\documentclass[a4paper,12pt]{article}
\usepackage{amssymb}
\usepackage{amsthm}
\theoremstyle{plain}
\newtheorem{Lemma}{Lemma}
\newtheorem{Proposition}{Proposition}
\newtheorem{Theorem}{Theorem}
\newtheorem{Corollary}{Corollary}
\theoremstyle{definition}

\newtheorem{Remark}{Remark}

\newcommand{\K}{\mathbb{K}}

\title{On the surjectivity of the canonical \\ Gaussian map 
for multiple coverings}
\author{Edoardo Ballico - Claudio Fontanari}
\date{}

\begin{document}
\maketitle
\abstract Here we investigate the canonical Gaussian map for higher 
multiple coverings of curves, the case of double coverings being 
completely understood thanks to previous work by Duflot. 
In particular, we prove that every smooth curve can be covered 
with degree not too high by a smooth curve having arbitrarily 
large genus and surjective canonical Gaussian map.
As a consequence, we recover an asymptotic version of a recent result 
by Ciliberto and Lopez and we address the Arbarello subvariety in the 
moduli space of curves.

\vspace{0.5cm}

\noindent
\textsc{A.M.S. Math. Subject Classification (2000)}: 14H47, 14H51, 14H55.


\section{Introduction}
Let $X$ be a smooth projective algebraic variety $X$ and $L$ be a line 
bundle on $X$. The \emph{Gaussian map} 
$$
\gamma_X(L): \bigwedge^2 H^0(X,L) \longrightarrow H^0(X, \Omega^1_X \otimes 
L^{\otimes 2})
$$  
is locally defined by the formula $\gamma_X(L)(\sigma \wedge \tau)=
\sigma d \tau - \tau d \sigma$ (see for instance \cite {Wahl:92}, p. 304). 
If $X=C$ is a curve and $L=K$ is the canonical line bundle on $C$ we say 
that $\gamma_C:=\gamma_C(K)$ is the \emph{canonical} Gaussian map for $C$. 
Even though many authors refer to it as to the \emph{Wahl map} after
Jonathan Wahl, the first steps towards the study of low gonality curves from 
the point of view of Gaussian maps were moved by the old Italian geometers: 
in particular, the rank of the so-called Wahl map for hyperelliptic curves 
was implicitly established by Beniamino Segre in his beautiful memoir 
\cite{Segre:27}, which goes back to 1927. The general case of double 
coverings was addressed by Jeanne Duflot in \cite{Duflot:94}. 
Gaussian maps for trigonal curves were 
studied independently by Ciro Ciliberto and Rick Miranda (see \cite
{CilMir:92}) and by James N. Brawner (see \cite{Brawner:95}), who turned 
his attention also to tetragonal curves (see \cite{Brawner:97}). 
The case of higher gonality has been recently settled by Ciliberto  
and Angelo Felice Lopez (see \cite{CilLop:02}). Here we recover an 
asymptotic version of their Theorem~(1.2) from the following more 
general result:  

\begin{Theorem}\label{main}
Let $b \ge 0$ be an integer and set $c := \max \{2b+13, 4b+3 \}$. 
For every smooth projective curve $B$ of genus $b$ and for every 
integer $g \ge 2 c^2 - 5 c + 7$ there is a smooth projective curve 
$C$ of genus $g$ with a morphism $C \to B$ of degree $d \le 2(c+1-2b)$ 
such that $\gamma_C$ is surjective. 
\end{Theorem}

Notice that Theorem~\ref{main} contrasts the hope expressed in 
\cite{CilLop:02} that the general $d:1$ cover of a general curve 
of genus $b \ge 1$ could be nonextendable with nonsurjective Wahl map. 

A remarkable advantage of our approach, which rests upon a powerful 
construction by Wahl, is that it applies also to coverings with prescribed 
ramifications. In particular, we obtain that a general curve 
in the subvariety $\overline{W}_{n,g}$ defined by Enrico Arbarello 
in \cite{Arbarello:74} and \cite{Arbarello:78} has surjective canonical 
Gaussian map (see Proposition~\ref{Arbarello} for a more precise statement).

We work over an algebraically closed field $\K$ of characteristic zero.

This research is part of the T.A.S.C.A. project of I.N.d.A.M., supported 
by P.A.T. (Trento) and M.I.U.R. (Italy).

\section{The results}
As mentioned above, the cornerstone of this work is a construction 
due to Wahl. The following statement 
is essentially \cite{Wahl:90}, Theorem 4.11; our condition (a) is slightly 
weaker than the original one, but the same proof works without any change.

\begin{Lemma}\label{Wahl} \textbf{(Wahl)} Let $C_i$ be a smooth projective 
curve of genus $g_i$ ($i=1,2$), $K_i$ a canonical divisor on $C_i$ and 
$D_i$ a divisor on $C_i$ of degree $d_i$ such that: 
\newline
(a) the general curve $C$ in the linear series $\vert p_1^*D_1+p_2^*D_2 \vert$ 
on $C_1 \times C_2$ is smooth and connected; \newline
(b) $d_i > \max(0,4-4g_i)$; \newline
(c) $K_i+D_i$ is normally generated on $C_i$ and $\gamma_{C_i}(\mathbf
{\mathcal O}_{C_i}(K_i+D_i))$ is onto; \newline
(d) either $g_1 \ge 2$ or $g_2 \ge 2$. \newline
Then $\gamma_C$ is surjective and 
$$
g(C)=d_1(g_1+d_2-1)+d_2(g_2-1)+1.
$$
\end{Lemma}

\begin{Remark}\label{numerical}
If $g_2 \ge 2$, $d_i \ge \max\{g_i+4,5 \}$ if $g_i \le 3$ and 
$d_i \ge 2g_i+1$ otherwise, then all hypotheses of Theorem~\ref{Wahl}
are satisfied. Indeed, (a) follows from Bertini's Theorem since both 
the $D_i$'s are very ample; $K_i+D_i$ is normally generated by 
\cite{Mumford:70} p.~55, Corollary to Theorem~6, since $d_i \ge 2g_i+1$; 
$\gamma_{C_i}(\mathbf{\mathcal O}_{C_i}(K_i+D_i))$ is surjective by 
\cite{B-E-L:91}, Theorem~1~(i), since $2g_i-2+d_i \ge 3g_i+2$.
\end{Remark}

\emph{Proof of Theorem~\ref{main}.} For every integer $g_1 \ge 2b$ there 
exists a smooth projective curve $C_1$ of genus $g_1$ with a morphism 
$C_1 \to B$ of degree $2$. We implement the construction of 
Theorem~\ref{Wahl} with such a $C_1$  and $g_1+d_2-1=c$. 
For $g_2=2$ we get $g(C)=d_1 c + d_2 + 1$ with $6 \le d_2 \le c+1-2b$ and 
$d_1 \ge 2c-5$, while for $g_2=3$ we obtain $g(C)=d_1 c + 2d_2 + 1$ with 
$7 \le d_2 \le c+1-2b$ and $d_1 \ge 2c-7$.
This way we cover the cases $g(C) = d_1c+r$ with $7 \le r \le c+6$, 
$d_1 \ge 2c-5$. Hence for every integer $g \ge 2 c^2 - 5 c + 7$ 
we have produced a smooth curve $C$ with a morphism of degree 
$2 d_2 \le 2 (c + 1 - 2 b)$ onto $B$ and from Lemma~\ref{Wahl} and 
Remark~\ref{numerical} it follows that $\gamma_C$ is surjective. 

\qed

\begin{Corollary} \textbf{(Ciliberto--Lopez)} For all integers $g \ge 280$ 
and $k \ge 28$ the general $k$-gonal curve of genus $g$ has surjective 
canonical Gaussian map.
\end{Corollary}

\begin{Remark}
We stress that our numerical assumptions are far away from being sharp 
(cf. \cite{CilLop:02}, Theorem~(1.2)). 
\end{Remark}

Let now $\mathcal{M}_g$ be the moduli space of curves of genus $g \ge 2$. 
As in \cite{Arbarello:74}, let $\overline{W}_{n,g}$ denote the closure 
in $\mathcal{M}_g$ of the locus of curves admitting a degree $n$ morphism 
onto $\mathbb{P}^1$ with a point of total ramification; in other words, 
$\overline{W}_{n,g}$ corresponds to smooth curves $C$ of genus $g$ 
possessing a Weierstrass point $p$ such that 
$\dim H^0(C, \mathcal{O}_C(np))\ge 2$. 

\begin{Lemma}\label{nonempty} 
Notation as in Theorem~\ref{Wahl}. 
Choose $q_2 \in C_2$ and let $D_2 := d_2 q_2$. If $C_1$ is hyperelliptic 
and $d_1 \ge 2 g_1$, then there exists $C \in \vert p_1^*D_1+p_2^*D_2 \vert$ 
such that $[C] \in \overline{W}_{2d_2,g}$. 
\end{Lemma}

\proof Take $q_1 \in C_1$ such that $g^1_2 = \vert 2 q_1 \vert$, 
let $q:= (q_1, q_2) \in C_1 \times C_2$ and consider the 
zero-dimensional scheme $Z := (q_1, d_2 q_2)$ on $\{ q_1 \} \times 
C_2$. We define the following three open subsets of 
$\mathcal{L} := \vert p_1^*D_1+p_2^*D_2(-Z)\vert$: 
\begin{eqnarray*}
U_1 &=& \{ D \in \mathcal{L}: D \cap \{ q_1 \} \times C_2 = Z \}\\
U_2 &=& \{ D \in \mathcal{L}: D \textrm{ is smooth outside of $q$} \}\\
U_3 &=& \{ D \in \mathcal{L}: D \textrm{ is smooth in $q$} \}\\
\end{eqnarray*}
If we show that $U_i \ne \emptyset$ for every $i = 1, 2, 3$, then 
$U_1 \cap U_2 \cap U_3 \ne \emptyset$, hence our claim follows. 

In order to check that $U_1 \ne \emptyset$, consider the natural 
exact sequence on $C_1 \times C_2$: 
\begin{eqnarray*}
0 &\to& \mathcal{O}(p_1^*D_1+p_2^*D_2 - \{q_1\}\times C_2) \to 
\mathcal{O}(p_1^*D_1+p_2^*D_2)\\ 
&\to& \mathcal{O}(p_1^*D_1+p_2^*D_2)\vert_{\{q_1\}\times C_2} \to 0. 
\end{eqnarray*}
Since $d_1 \ge 2g_1$, we have $H^1(C_1, \mathcal{O}(D_1-q_1))=0$ and 
$H^1(C_1 \times C_2, \mathcal{O}(p_1^*D_1+p_2^*D_2 - \{q_1\}\times C_2)=0$ 
by the K\"unneth formula. Hence 
$$
H^0(C_1 \times C_2,\mathcal{O}(p_1^*D_1+p_2^*D_2)) \to 
H^0(\{q_1\}\times C_2,\mathcal{O}(p_1^*D_1+p_2^*D_2)) \to 0  
$$
and if we take a section of $\mathcal{O}(D_2) = \mathcal{O}(d_2 q_2)$ 
vanishing on $d_2 q_2$, any section in its non-empty counterimage 
has an element of $U_1$ as its zero divisor.

Next, we have $U_2 \ne \emptyset$ by Bertini's Theorem; finally, it 
is easy to produce a divisor in $\mathcal{L}$ which is smooth in $q$: 
just take $\{q_1\} \times C_2 + D$ with $D \in p_1^*D_1(-q_1)+p_2^*D_2$
not passing through $q$. 

\qed

We recall from \cite{Arbarello:78} that $\overline{W}_{n,g}$ is an 
irreducible subvariety of $\mathcal{M}_g$, so it makes sense to consider 
its general element. 

\begin{Proposition}\label{Arbarello}  
For all integers $g \ge 280$ and $28 \le n \le g$, a general curve in 
$\overline{W}_{n,g}$ has surjective canonical Gaussian map.
\end{Proposition}

\proof By \cite{Arbarello:74}, Theorem~(3.18), we have 
$\overline{W}_{n-1,g} \subset \overline{W}_{n,g}$ 
for every $2 \le n \le g$. Hence from the proof of Theorem~\ref{main} 
and Lemma~\ref{nonempty} it follows that there exists 
$C \in \overline{W}_{n,g}$ with $\gamma_C$ surjective. 
Now the claim follows from the semicontinuity of the corank of the 
canonical Gaussian map.

\qed

\vspace{0.5cm}
\noindent Edoardo Ballico \newline
Universit\`a degli Studi di Trento \newline
Dipartimento di Matematica \newline
Via Sommarive 14 \newline
38050 POVO (Trento) \newline
Italy \newline
e-mail: ballico@science.unitn.it

\vspace{0.5cm}
\noindent
Claudio Fontanari \newline
Universit\`a degli Studi di Trento \newline
Dipartimento di Matematica \newline
Via Sommarive 14 \newline
38050 POVO (Trento) \newline
Italy \newline
e-mail: fontanar@science.unitn.it

\end{document}